\documentclass[reqno]{amsart}
\usepackage{amssymb}
\usepackage{amsmath}

\makeatletter
\@addtoreset{equation}{section}
\makeatother

\renewcommand\thefigure{\thesection.\@arabic\c@figure}
\renewcommand\thetable{\thesection.\@arabic\c@table}

\newtheorem{theorem}{Theorem}[section]

\newtheorem{proposition}[theorem]{Proposition}

\newtheorem{remark}[theorem]{Remark}

\def\X{D}
\def\d{\delta}
\def\qed{\hfill$\nabla$}

\usepackage{ulem}
\usepackage{color}
\definecolor{nb}{rgb}{.6,.176,1}
\definecolor{sienna}{rgb}{.52,.222,.176}
\definecolor{darkgreen}{rgb}{0,.5,0}

\begin{document}

\author{Siva Athreya, Sunder Sethuraman, and Balint Toth}

\address{Siva Athreya \\ Indian Statistical Institute
8th Mile Mysore Road  \\ Bangalore 560059, India\\
e-mail: \rm \texttt{athreya@isibang.ac.in}}

\thanks{Research supported in part by a CSIR Grant and Homi Bhaba Fellowship.}

\address{Sunder Sethuraman\\
396 Carver Hall\\
Department of Mathematics\\
Iowa State University\\
Ames, IA \ 50011\\
e-mail:  \rm \texttt{sethuram@iastate.edu} }

\thanks{Research supported in part by NSA-Q H982301010180,
NSF-DMS 0906713.}

\address{Balint Toth \\ Institute of Mathematics\\
Technical University Budapest\\
Egry Jozsef u. 1\\
H-1111 Budapest, Hungary\\
e-mail: \rm \texttt{balint@math.bme.hu} }

\title[]
{On the range, local times and periodicity of random walk on an interval}

\begin{abstract}
 The range, local times, and periodicity of symmetric, weakly asymmetric and asymmetric random walks at the time of exit from a strip with $N$ locations are considered.  Several results on asymptotic distributions are obtained.
\end{abstract}



\maketitle

\section{Introduction}


In this article, we study the range, local times,
and periodicity or ``parity'' statistics of nearest-neighbor
symmetric, weakly asymmetric, and asymmetric random walks up to the
time of exit from an interval of $N$ sites.
We derive several associated scaling limits which appear curious, some which
connect with the entropy of an exit distribution, generalized
Ray-Knight constructions, and Bessel and Ornstein-Uhlenbeck square processes, among other
objects.

The study of the range of random walk is of course an old subject.
However, examining the range and related statistics at the time the random walk leaves an interval, although a simple, natural concern, seems unexplored.  We refer to Bass-Chen-Rosen \cite{BCR}[Ch. 2] and Den Hollander-Weiss \cite{DenHollander_Weiss} for exhaustive references on the range and related statistics of random walk in various settings.

From another view, indeed our initial motivation for this problem, the study of the range and other structures of random walk when it exits an interval can be thought of as a stochastic version of the ``locker'' problem, popular in university curriculum:
Suppose there is a hallway of lockers labeled from $1$ to $N$, for $N\geq
1$, which are initially closed.  Let persons $L$, for $L\geq 1$, walk through the hallway, toggling every $L$th locker, that is opening it if closed and closing it if open.
The question is then to find out those lockers which will be open after the first $N$ people walk through.  The lockers whose labels are the squares, $1$, $4$, $9$, etc., are exactly those with an odd number of factors.  Consequently, these lockers are 
 the open ones.   Other variations of this problem can be found in
Tanton \cite{tanton} and references therein.


In our random walk setting, we can imagine each site in the interval to be either open or closed, and the random walker toggling a site on each visit (from open to closed and vice versa) before it exits the interval.  In comparison to the  ``locker'' problem, we address  the following questions:-

\begin{itemize}
\item[(1)] What fraction of sites will be visited when the walker exits, e.g. the range?
\item[(2)] How many times will each site be visited before exit, e.g local times across sites?
\item[(3)] And, given a set of sites 
that have been visited, what is the joint distribution of their open status at the time of exit, e.g parity of the visits to points in the interval?
\end{itemize}

 The specific answers naturally depend on the type of random walk
 considered.  A goal of the paper is to see how the behaviors under
 symmetric and asymmetric walks are interpolated in terms of weakly
 asymmetric walks.

For the first question, we derive 
the limiting distribution for the range
(Proposition \ref{rangeprop}), and observe as a consequence,
which seems surprising, that the scaled range, when starting at
random, is uniformly distributed on $[0,1]$ no matter the dynamics
(Proposition \ref{uniformprop}).  Also, curious values for the expected scaled range
under symmetric walks, and the chance a given point is in the range,
when starting at random are found (Remarks
\ref{entropy} and \ref{point_visitedrmk}).

For the second question, we find the scaling limit of the local times
through a ``Ray-Knight'' construction involving Bessel and
Orstein-Uhlenbeck squared processes (Propositions \ref{Besqprop},
\ref{OUprop} and \ref{asym_limprop}).

For the third question, we show that
the parities of well-separated points, given that they are visited, are independent and
identically distributed Bernoulli variables, and fair in the
symmetric/weakly asymmetric case, and biased in the asymmetric
situation (Proposition \ref{thm1_question3} and \ref{thm2_question3}).

\vskip .1cm
{\bf Set up:-}
Let $\mathcal T_{N} = \{0,1,2,\ldots, N\}.$  Let $X_n$ be the position
of a random walk on $\mathcal T_N$ at times $n\geq 1$. At each
time step, the walk moves to the nearest point to its left
(right) with probability $q_N$ ($p_N$) where $p_N+q_N=1$.  The walk stops the moment it
is at either $0$ or $N$. When $p_N=q_N=1/2$, the walk is of course referred to as the symmetric random
walk.  When $q_N = 1/2 -c/N$ (and so $p_N = 1/2+c/N$) for some
constant $c>0$ and $N$ large enough so that $0< p_N,q_N<1$, we say the
walk is weakly asymmetric.  When $q_{N}= q < p = p_{N} $, the walk
is asymmetric.

Define $T_a = \inf \{n \geq 1: X_{n} = a\}$ as the hitting time of $a\in \mathcal T_N$.  Then, $\tau_N = T_0\wedge T_N$ is the ``exit'' time from the strip.
Clearly, starting from $1\leq x\leq N-1$, $\tau_N$ is
finite: $P_x(\tau_N<\infty)=1$ where we denote
$P_x(A)=P(A|X_0=x)$ as the conditional probability of the event $A$
with respect to the walk starting from $X_0=x$.

Then, the number of visits
to $y \in \mathcal T_{N}$ before exiting is
$G(y)  =  \sum_{k=0}^{\tau_N} 1_{y}(X_k)$.
Hence, the event $y$ is visited at all corresponds to $G(y)\geq 1$.
In this case, we say the parity of $y$ is ``even'' (locker $y$ is closed) if $G(y)\geq 1$ and $G(y)=0 \ {\rm mod}_2$.  Correspondingly, the parity of $y$ is ``odd'' (locker $y$ is open) when $G(y)\geq 1$ and
                $G(y)=1\ {\rm mod}_2$.\vskip .1cm

The plan of the article is to address questions (1),(2) and (3) in
sections \ref{distributional}, \ref{ray-knight}, and \ref{independent}
respectively.

\section{Question 1: Range of random walk in $\mathcal T_N$}
\label{distributional}

In this section, we obtain
distributional limits of the range up to the exit time when starting
from a point, and at random in
subsections \ref{distributional1}, \ref{random}.

\subsection{The range starting from a point}
\label{distributional1}

Denote $R_N$ as the number of locations visited before exit, the range of the walk on $\mathcal T_N$, that is
$$R_N \ = \ \#\{y\in \mathcal T_N: G(y)\geq 1\}.$$
Observe, when starting from $[ \alpha N]$, necessarily $[ \alpha N\wedge (1-\alpha) N] \leq R_N\leq N$.


\begin{proposition}
\label{rangeprop} Let $X_{0}= [ \alpha N]$ for $0<\alpha<1$.   For symmetric and weakly asymmetric walks, $R_N/N$ converges in distribution to absolutely continuous measures on $[0,1]$, respectively $G_{0,\alpha}$ and $G_{c,\alpha}$ defined in (\ref{G_symmetric}) and (\ref{G_weakly}). For asymmetric walks, $R_N - [(1-\alpha)N] \Rightarrow Z$ where $Z$ is Geometric$(q/p)$.

 \end{proposition}

\proof First, we observe, for $0<\beta<1$, starting from location $x=[\alpha N]$, since the motion is nearest-neighbor,
\begin{eqnarray*}
\{R_N \geq \beta N\} & = & \{R_N \geq \beta N, \tau_N = T_0\} \cup \{R_N \geq \beta N, \tau_N = T_N\}\\
&=& \{T_{[ \beta N]} < T_0<T_N\} \cup \{T_{N-[ \beta N]}<T_N<T_0\}.
\end{eqnarray*}
We now specialize arguments to the three types of random walks.
\vskip .1cm

{\it Symmetric Walk:-}
When the walk is symmetric $p_N=q_N=1/2$, recall the standard Gambler's ruin identity: For $a,b,z \in \mathcal T_{N},$ such that $a<z<b  $,
\begin{equation}\label{gamblers_symmetric}
P_z(T_a<T_b) \ = \ \frac{b-z}{b-a}.\end{equation}

For $\beta> \alpha$, compute
 \begin{eqnarray*}
P_{[ \alpha
N]}(R_N \geq [ \beta N], \tau_N = T_0)
&=& \frac{[ \alpha
N]}{[ \beta N]}\frac{N-[ \beta N]}{N} \
\rightarrow\ \frac{\alpha(1-\beta)}{\beta}.\end{eqnarray*}

When, $\beta> 1-\alpha$, we have
\begin{eqnarray*} P_{[ \alpha
N]}(R_N \geq [ \beta N], \tau_N = T_N)
&=& \frac{N-[ \alpha
N]}{[ \beta N]}\frac{N-[ \beta
N]}{N}\ \rightarrow \
\frac{(1-\alpha)(1-\beta)}{\beta}.\end{eqnarray*}

Putting these expressions together, along with simple calculations, we have
$$\lim_{N\uparrow \infty}P_{[ \alpha
N]}(R_N/N\geq \beta) \ = \ \left\{\begin{array}{rl}
1&\ {\rm when \ }0\leq \beta\leq \alpha \wedge (1-\alpha)\\
\frac{\alpha \wedge (1-\alpha)}{\beta} &\ {\rm when \
}\alpha \wedge (1-\alpha) <\beta
<\alpha \vee (1-\alpha)\\
\frac{1-\beta}{\beta}&\ {\rm when \ }
\alpha \vee (1-\alpha) \leq \beta \leq 1\\
0&\ {\rm when \ }\beta>1.\end{array}\right.$$
The right-side defines a distribution $G_{0,\alpha}$, supported on $[\alpha\wedge (1-\alpha), 1]$ whose density
\begin{equation}
\label{G_symmetric}
g_{\alpha}(\beta) \ = \ \left\{\begin{array}{rl}
\frac{\alpha \wedge (1-\alpha)}{\beta^2}& \ {\rm for \
}\alpha \wedge (1-\alpha) <\beta
<\alpha \vee (1-\alpha)\\
\frac{1}{\beta^2}& \ {\rm for \ }\alpha \vee (1-\alpha)\leq \beta \leq 1\\
0&\ {\rm otherwise.}\end{array}\right.\end{equation}

%
\vskip .1cm

{\it Weakly-asymmetric walk:-} In the weakly asymmetric case, $q_N= 1/2 - c/N$ and $p_N = 1/2 + c/N$ with $c>0$, let
$$s_N \ := \ \frac{q_N}{p_N} \ = \ \frac{1/2 -c/N}{1/2+c/N} \ = \ 1-\frac{4c}{N} + O(N^{-2}).$$
  The corresponding gambler's ruin identity becomes, for $a<z<b$,
\begin{equation}\label{gamblers_weakly}
P_z(T_a<T_b) \ = \ \frac{(q_N/p_N)^z - (q_N/p_N)^b}{(q_N/p_N)^a - (q_N/p_N)^b}.\end{equation}


Then, following the symmetric argument, when $\beta> \alpha$,
\begin{eqnarray*}
P_{[ \alpha
N]}(R_N \geq [\beta N], T_0= \tau_N)
&=& \frac{1-s_N^{[ \alpha N]}}{1-s_N^{[ \beta N]}} \frac{s_N^{[ \beta N]} - s_N^N}{1-s_N^N}\\
&\rightarrow& \frac{1-e^{-4\alpha c}}{1-e^{-4\beta c}}\frac{e^{-4\beta c}-e^{-4c}}{1-e^{-4c}}\ := \ A_1(\alpha,\beta,c).
\end{eqnarray*}
When $\beta>1-\alpha$,
\begin{eqnarray*}
P_{[ \alpha
N]}(R_N\geq [\beta N], T_N=\tau_N)
&=& \frac{s_N^{[ \alpha N]}-s_N^N}{s_N^{N-[ \beta N]}-s_N^N}\frac{1-s_N^{N-[ \beta
N]}}{1-s_N^N}\\
&\rightarrow &
\frac{e^{-4\alpha c} - e^{-4c}}{e^{-4(1-\beta)c}-e^{-4c}}\frac{1-e^{-4(1-\beta) c}}{1-e^{-4c}} \ := \ A_2(\alpha,\beta,c).\end{eqnarray*}

Noting
$$
\lim_{N\uparrow\infty}P_{[ \alpha N]}(T_N<T_0) \ =\  \frac{1-e^{-4\alpha c}}{1-e^{-4c}},\ \ {\rm and \ \ }
A_1+A_2 \ = \ \frac{e^{-4\alpha c}(e^{4c(1-\beta)}-1)}{1-e^{-4c\beta}},
$$
we have
\begin{equation}
\label{G_weakly}
\lim_{N\uparrow \infty}P_{[ \alpha
N]}(R_N/N\geq \beta) \ = \ \left\{\begin{array}{rl}
1&\ {\rm when \ }0\leq \beta\leq \alpha \wedge (1-\alpha)\\
\frac{1-e^{-4c\alpha}}{1-e^{-4c\beta}} &\ {\rm when \
}\alpha  \leq \beta
\leq 1-\alpha\\
\frac{1-e^{4c(1-\alpha)}}{1-e^{4c\beta}}& \ {\rm when \ } 1-\alpha \leq \beta \leq \alpha\\
\frac{e^{-4\alpha c}(e^{4c(1-\beta)}-1)}{1-e^{-4c\beta}}&\ {\rm when \ }
\alpha \vee (1-\alpha) \leq \beta \leq 1\\
0&\ {\rm when \ }\beta>1.\end{array}\right.\end{equation}
The right-side defines a distribution $G_{c,\alpha}$, supported on $[\alpha\wedge (1-\alpha), 1]$, whose density, although messy, can be easily derived.

\vskip .1cm
{\it Asymmetric walk:-} In the asymmetric case, $q_N=q$, $p_N=p$, and $p>q$, and it is not difficult to see that we cannot go left too many times.  The gambler's ruin identity (\ref{gamblers_weakly}) also holds in this case and,
for $x= [\alpha N]$, $P_{x}(T_0<T_N) = \exp\{-CN\}$ for some constant $C>0$.

To complete the proof, for integers $z\geq 0$, compute
\begin{eqnarray*}
P_x(R_N \geq N-x + z) &=& P_x(T_N<T_0, R_N\geq N-x +z) + o(1)\\
&=& P_x(T_{x-z}<T_N<T_0) + o(1)\\
&=& (q/p)^{z} + o(1).
\end{eqnarray*}
\qed

\begin{remark}
\label{entropy}
 \rm

 1.  As expected, $G_{c,\alpha}$ interpolates between the symmetric and asymmetric cases:  Namely, as $c\downarrow 0$, $G_{c,\alpha}\Rightarrow G_{0,\alpha}$, and as $c\uparrow\infty$, $G_{c,\alpha}$ converges to the constant $1-\alpha$.

2.  It is curious to observe, for symmetric walks, that starting from $x = [ \alpha N]$, with $\alpha \in (0,1/2]$, the expected range
$$
\int_0^1 \beta g_\alpha(\beta)d\beta
\ = \ \int_{\alpha}^{1-\alpha}  \beta\frac{\alpha}{\beta^2}d\beta + \int_{1-\alpha}^1 \beta\frac{1}{\beta^2}d\beta
\ =\  -(1-\alpha)\log(1-\alpha)
-\alpha\log(\alpha)
$$
        is the entropy of the exit distribution $\langle 1-\alpha,\alpha\rangle$ where $1-\alpha$ is the probability of leaving by the left endpoint, and $\alpha$ the chance of exiting right!  The maximum value $\log 2$ occurs when $\alpha=1/2$.


3.  For symmetric and weakly asymmetric walks, the limit distributions may also be derived in terms of Brownian motion and diffusion estimates.
\end{remark}

\subsection{The range when starting at random} \label{random} We derive now the limiting law of $R_N/N$ when starting at random, that is the uniform distribution on $\mathcal T_N$.  It seems nonintuitive that the limit law is U$[0,1]$ no matter the type of random walk.

\begin{proposition}
\label{uniformprop}
For symmetric, weakly asymmetric and asymmetric random walk, when starting at random in $\mathcal T_N$, $R_N/N$ converges weakly to the uniform distribution $U[0,1]$.
\end{proposition}

\proof  Suppose our starting point was random. In the symmetric and weakly asymmetric cases, the limiting distribution of ${R_{N}}/{N}$, from straightforward considerations, is found by integrating the density $g_\alpha$ and tail of $G_{c,\alpha}$ with respect to $\alpha$ (denoted by $G_{c,\alpha}([\beta,1])$).

In the symmetric case, when $\beta\leq 1/2$,
\begin{eqnarray*}
\int_0^1 g_\alpha(\beta)d\alpha &=& \int_0^\beta \frac{\alpha}{\beta^2} d\alpha + \int_\beta^{1-\beta} 0\; d\alpha + \int_{1-\beta}^1 \frac{1-\alpha}{\beta^2} d\alpha \ = \ 1.
\end{eqnarray*}
But, also, when $\beta>1/2$,
\begin{eqnarray*}
\int_0^1 g_\alpha(\beta)d\alpha &=& \int_0^{1-\beta} \frac{\alpha}{\beta^2}d\alpha + \int_{\beta}^{1-\beta} \frac{1}{\beta^2}d\alpha + \int_{\beta}^1 \frac{1-\alpha}{\beta^2}d\alpha \ = \ 1.
\end{eqnarray*}

On the other hand, in the weakly asymmetric case, we have, when $\beta\leq 1/2$,
\begin{eqnarray*}\int_0^1 G_{c,\alpha}([\beta,1])d\alpha &=&
\int_0^\beta \frac{1-e^{-4c\alpha}}{1-e^{-4c\beta}}d\alpha + \int_\beta^{1-\beta}1\; d\alpha \\
&&+ \int_{1-\beta}^1 \frac{e^{-4c\beta}(e^{4c(1-\alpha)} - 1)}{e^{-4c(1-\beta)}-e^{-4c}}d\alpha\ = \ 1-\beta.
\end{eqnarray*}
Similarly, when $\beta>1/2$, $\int_0^1 G_{c,\alpha}([\beta,1])d\alpha$ equals
\begin{eqnarray*}
&&\frac{1}{1-e^{-4c\beta}}\Big[\int_0^{1-\beta}{1-e^{-4c\alpha}}d\alpha + \int_{1-\beta}^\beta e^{-4c\alpha}(e^{4c(1-\beta)}-1)d\alpha\\
 &&\ \ \ \ \ \ \ \ \ \ \ \ \ \ \ \ \ \ \ \ \ \ \ \ \ \ \ \ \  + \int_\beta^1 e^{-4c\beta}(e^{4c(1-\alpha)}-1)d\alpha\Big] \ = \ 1-\beta.
\end{eqnarray*}
Consequently, for symmetric and weakly asymmetric walks, the limiting distribution is $U[0,1]$ when the starting position is uniformly chosen.

However, in the asymmetric case, from Proposition \ref{rangeprop}, $R_N/N \rightarrow 1-\alpha$ in probability starting from $x=[ \alpha N]$.  Then, starting at random in $\mathcal T_N$, we have that $R_N/N\rightarrow Y$ in probability where $Y$ is a $U[0,1]$ distributed random variable. \qed

\begin{remark}
\label{point_visitedrmk}
\rm
One might ask, on the other hand, with what probability a point $y= [ \beta N]$ belongs to the range when starting at random.  This is the same as asking when $y$ is visited by the walk.  For symmetric walk, it is not difficult to use the gambler's ruin identity (\ref{gamblers_symmetric}) to see,
as $N\uparrow\infty$, that the probability tends to
$$\int_\beta^1 \frac{1-\alpha}{1-\beta}d\alpha +
\int_0^\beta \frac{\alpha}{\beta}d\alpha\ = \ \frac{1}{2}.$$
It seems curious that the limit does not depend on $\beta$.

For asymmetric walk, starting from $[ \alpha N]$, when $\alpha>\beta$, the point $y$ cannot be reached with positive probability in the limit.  Then, the chance $y$ belongs to the range, when starting at random, is $\beta$.

For weakly asymmetric walks, using (\ref{gamblers_weakly}), the limit is $\frac{\beta}{1-e^{-4c\beta}} - \frac{1-\beta}{1-e^{4c(1-\beta)}}$ which interpolates between the other cases as $c\downarrow 0$ and $c\uparrow\infty$
\end{remark}

\section{Question 2:  Characterization of local times}
\label{ray-knight}
To capture the local times of the random walk before its exit, we use the ``Ray-Knight'' or ``Kesten-Kozlov-Spitzer'' representation, and some martingale characterizations.  Our treatment and proofs will be similar to those in Toth \cite{Toth} which considered certain self-interacting random walks.

Let $0 < \alpha < 1$.  Suppose the walk starts at $[\alpha N]$, and exits at the right endpoint $N$. Let $\zeta^N_{j}$ be the number of left crossings of the bond $(N-(j-1),N-j)$ before exit.  Then, $\zeta^N_0=0$, and $\zeta^N_1$ is distributed as $D_N$, a Geometric$(q_N)$ random variable minus $1$,
$P(D_N=n) = p_Nq_N^{n}$ for $n\geq 0$.  In the following, we drop the script $N$.


Let $\{\xi_{j,i}\}_{i,j\geq 0}$ be i.i.d. random variables with distribution $D_{N}$.  A moment's thought convinces that $\{\zeta_j\}_{0\leq j\leq N}$ is a Markov chain with representation
\begin{equation}
\label{representation}\zeta_{j+1} \ = \ \left\{\begin{array}{rl}\sum_{i=0}^{\zeta_j} \xi_{j,i}& \ {\rm for \ } 0\leq j<[(1-\alpha)N]\\
\sum_{i=1}^{\zeta_j} \xi_{j,i} & \ {\rm for \ } [(1-\alpha)N]\leq j\leq N-1\end{array}\right.
\end{equation}
such that
\begin{equation}
\label{restriction}\zeta_j=0 {\rm \ \ \ for\  some\ \ \ } [(1-\alpha)N]\leq j<N,\end{equation}
with the convention that empty sums vanish.

Note that for $j< [(1-\alpha)N]$, the sum starts with index $i=0$ since, even if $\zeta_j=0$, given exit at the right, the walk must visit locations $[\alpha N]\leq x\leq N$ and may have left crossings of $(x-1,x)$.
However, for $j\geq [(1-\alpha)N]$, since the walk is not guaranteed to visit sites to the left of $[\alpha N]$, $\zeta_j$ is the size of a Branching process, with initial value $\zeta_{[(1-\alpha)N]}$, which must vanish before time $j=N$.

Then, the local time of the walk is $$G(y) \ =\  \left\{\begin{array}{rl}
\zeta_{N-y}& \ {\rm for \ }0 \leq y < [\alpha N]\\
\zeta_{N-y} +1 & \ {\rm for \ }[\alpha N] \leq y\leq N. \end{array}\right.
$$

In the following, to analyze $\{\zeta_j\}_{0\leq j\leq N}$, it will be helpful to consider the Markov chain $\eta_j$, such that $\eta_0=0$ and $\eta_1\stackrel{d}{=}D_N$, for which representation (\ref{representation}) holds in terms of the variables $\{\xi_{j,i}\}_{i,j\geq 0}$, but {\it without} the restriction (\ref{restriction}).

When the walk exits at the left endpoint $0$, one considers an analogous Markov chain $\tilde\zeta_j$, corresponding to right-crossings of $(j,j+1)$, where the representation and restriction are reversed.  Namely, let $\tilde D_N$ be a ${\rm Geometric}(p_N)$ random variable minus $1$, $P(\tilde D_N = n) = q_Np_N^{n}$ for $n\geq 0$.  Define $\tilde\zeta_0=0$, $\tilde \zeta_1 \stackrel{d}{=}\tilde D_N$, and
$$
\zeta_{j+1} \ = \ \left\{\begin{array}{rl}\sum_{i=0}^{\zeta_j} \xi_{j,i}& \ {\rm for \ } 0\leq j<[\alpha N]\\
\sum_{i=1}^{\zeta_j} \xi_{j,i} & \ {\rm for \ } [\alpha N]\leq j\leq N-1\end{array}\right.
$$
such that
$\zeta_j=0$ for some $[\alpha N]\leq j<N$.  The local time of the walk in this case is $G(y) = \tilde\zeta_{y}$ for $[\alpha N]< y\leq N$ and $G(y) = \tilde\zeta_y +1$ for $0\leq y\leq [\alpha N]$.
  Here also it will be of use to define analogously a Markov chain $\tilde\eta_j$ satisfying $\tilde\eta_0=0$, $\tilde\eta_1\stackrel{d}{=}\tilde D_N$, and the reversed representation but without the restriction that the chain must vanish for $[\alpha N]\leq j<N$.

Finally, define $Y_N(t) = \frac{1}{N}\eta_{[Nt]}$ and $\tilde Y_N(t) = \frac{1}{N}\tilde \eta_{[Nt]}$ for $0\leq t\leq 1$,
 and suppose that $Y_N(0)=\tilde Y_N(0)=0$.


\subsection{Symmetric walks}

Consider the following processes.  Let $Z_0=0$, and define
$$ Z_t \ =  \ \left\{\begin{array}{rl}
t + \int_{0}^{t}\sqrt{2Z_s}dB_s & \ {\rm for \ } 0\leq t\leq 1-\alpha\\
&\\
Z_{1-\alpha} + \int_{1-\alpha}^{t}\sqrt{2Z_s}dB_s & \ {\rm for \ } 1-\alpha \leq t\leq 1.
\end{array}\right.
$$
Observe that
$Z_t$ for $0\leq t\leq 1-\alpha$ is the same in law as Besq$^2(t/2)$ process, and a Besq$^0(t/2)$ process for $1-\alpha\leq t\leq 1$ (cf. Revuz-Yor \cite{RY} for more on the processes Besq$^\delta$ $dX_t = \delta dt + 2\sqrt{X_t}dB_t$).  Let $\tau^R_0$ be the first time $Z_t$ hits $0$ after time $t=1-\alpha$.  Note that $Z_t$ remains at value $0$ after time $\tau^R_0$.

Define also
$\tilde Z_t$ where $\tilde Z_0=0$ and
$$\tilde Z_t \ =   \ \left\{\begin{array}{rl}
t + \int_0^t\sqrt{2\tilde Z_s}dB_s & \ {\rm for \ } 0\leq t\leq \alpha\\
\tilde Z_\alpha + \int_\alpha^t\sqrt{2\tilde Z_s}dB_s & \ {\rm for \ } \alpha \leq t\leq 1.
\end{array}\right.
$$
Let also $\tau^L_0$ be the time $\tilde Z_t$ reaches $0$ after time $t=\alpha$.  Here, also, $\tilde Z_t \equiv 0$ for $t\geq \tau^L_0$.

It will turn out that $Z_t$ and $\tilde Z_t$ will be identified respectively, as the scaling limits of the local times when the random walk exits at the right and left endpoints of the interval.  The important point in this identification is the next result.


\begin{proposition}
\label{mainprop}
For symmetric walk starting from $x=[\alpha N]$, we have $$Y_{N}(t) \ \Rightarrow \ Z(t) \ \ \ {\rm and \ \ \ } \tilde Y_N(t) \ \Rightarrow \ \tilde Z(t)$$ in $D[0,1]$, in the sup topology.
\end{proposition}

 Instead of proving Proposition \ref{mainprop}, which follows steps as in Toth \cite{Toth}, we prove Proposition \ref{mainprop_OU} in the next subsection, with respect to weakly asymmetric random walks, dealing with squared Ornstein-Uhlenbeck processes which are less standard.

Now, with Proposition \ref{mainprop} in hand, since $Y_N(t)$ and $\tilde Y_N(t)$ converge respectively to $Z_t$ and $\tilde Z_t$ in the sup topology, it follows that the conditional distributions of $Y_N(t)$ given $\eta_j$ vanishes for $j\geq [(1-\alpha)N]$ and $\tilde Y_N(t)$ given $\tilde \eta_j$ vanishes for $j\geq [\alpha N]$ converge to the conditional distributions of $Z_t$ given that $1-\alpha\leq \tau^R_0<1$ and $\tilde Z_t$ given that $\alpha\leq \tau^L_0<1$.

Hence, from this discussion, the following characterization holds for the local times of the walk up to time of exit.  Recall that $1-\alpha$ and $\alpha$ are the exit probabilities of right and left exit respectively.
\begin{proposition}
\label{Besqprop} For symmetric walk starting from $[\alpha N]$, the local times
$$G([Nt])/N \ \Rightarrow \ \alpha \mu^R + (1-\alpha) \mu^L$$
 where $\mu^R$ is the law of the process $Z_{1-t}$ conditioned on $1-\alpha\leq \tau^R_0<1$, and $\mu^L$ is the law of the process $\tilde Z_t$ conditioned on $\alpha\leq \tau^L_0<1$.
\end{proposition}

\vskip .1cm

\subsection{Weakly asymmetric walks}

The development of the local time structure is similar to the symmetric case.  Corresponding to right exit, $E D_{N} = q_N/p_N =1-\frac{4c}{N+2c}$ and ${\rm Var}(\X_{N})=q^2_N/p^2_N + q_N/p_N = 2- \frac{12c}{N+2c} +\frac{16c^{2}}{(N+2c)^{2}}$.  Define the process $Z_t^c$ by
$Z^c_0=0$, and
$$
Z^c_t \ = \ \left\{\begin{array}{rl}
\int_{0}^{t}(1-4cZ^c_s)ds + \int_{0}^{t} \sqrt{2Z^c_s}dB_s & \ {\rm for \ } 0\leq t\leq 1-\alpha\\
&\\
Z^c_{1-\alpha}- \int_{1-\alpha}^{t}4cZ^c_s ds+ \int_{1-\alpha}^{t}\sqrt{2Z^c_s}dB_s & \ {\rm for \ } 1-\alpha \leq t\leq 1
\end{array}\right.
$$
Note $2(Z^c_t+1)$ and $2(Z^c_t -t)$ are the squares of the Ornstein-Uhlenbeck process $dX_t = -4cX_tdt +\sqrt{2}dB_t$ for $0\leq t\leq 1-\alpha$ and $1-\alpha\leq t\leq 1$ respectively.


Also, with respect to left exit, $E \tilde D_{N} = p_N/q_N = 1+4c/N + O(N^{-2})$ and ${\rm Var}(\tilde D_{N}) = 2+O(N^{-1})$.  Define
$\tilde Z^c_t$ by $\tilde Z^c_0=0$ and
$$\tilde Z^{c}_t \ = \ \left\{\begin{array}{rl}
\int_{0}^{t}(1+4c\tilde Z^c_s)ds + \int_{0}^{t}\sqrt{2\tilde Z^{c}_s}dB_s & \ {\rm for \ } 0\leq t\leq \alpha\\
\tilde Z^{c}_{\alpha} + \int_{\alpha}^{t}\sqrt{2\tilde Z^{c}_s}dB_s & \ {\rm for \ } \alpha \leq t\leq 1.
\end{array}\right.
$$

As before, let $\hat\tau^R_0$ be the first time after $t=1-\alpha$ that $Z^c_t$ reaches $0$, and $\hat\tau^L_0$ be the first time after $t=\alpha$ that $\tilde Z^c_t$ hits $0$.

Analogous to the symmetric random walk case, we show that $Z^c_t$ and $\tilde Z^c_t$ are the scaling limits of the local times when the weakly asymmetric random walk exits at the right and left endpoints respectively.

\begin{proposition}
\label{mainprop_OU}
For the weakly asymmetric random walk starting from $x=[\alpha N]$, we have $$Y_{N}(t) \ \Rightarrow \ Z^c(t) \ \ \ {\rm and \ \ \ } \tilde Y_N(t) \ \Rightarrow \ \tilde Z^c(t)$$ in $D[0,1]$, in the sup topology.
\end{proposition}

The same argument as in the symmetric case allows to deduce the the following characterization.

\begin{proposition}
\label{OUprop}
For the weakly asymmetric walk, starting from $x=[\alpha N]$, the local times satisfy
$$G([Nt])/N \ \Rightarrow \ R(\alpha) \mu_c^R + (1-R(\alpha)) \mu_c^L$$
 where $\mu_c^R$ is the law of the process $Z^c_{1-t}$ conditioned on $1-\alpha\leq \hat\tau^R_0<1$, and $\mu_c^L$ is the law of the process $\tilde Z^c_t$ conditioned on $\alpha\leq \hat\tau^L_0<1$.  Here, $R(\alpha) = (1-e^{-4c\alpha})/(1-e^{-4c})$ is the exit probability to the right.
 \end{proposition}

{\it Proof of Proposition \ref{mainprop_OU}.}  Here, we argue that $Y_N(t) \Rightarrow Z^c(t)$ which corresponds to ``left crossings.''  The argument for $\tilde Y_N(t)\Rightarrow \tilde Z^c(t)$ is similar.

The proof naturally separates into two parts corresponding to when $j\leq [(1-\alpha)N]$ and $j\geq [(1-\alpha)N]$.
The strategy will be to use martingale decompositions of the Markov chain $\{\eta_j\}_{j\geq 0}$.
Define, for $[Nt]\leq [(1-\alpha) N]$, the martingale and its quadratic variation,
\begin{eqnarray*}
M_N(t) & = & \eta_{[Nt]}-\eta_0 - \sum_{j=0}^{[Nt]-1} \left( E[\eta_{j+1}|\eta_j]-\eta_j\right)\\
\langle M_N(t) \rangle &=& \sum_{j=0}^{[Nt]-1} E\big[ \big(\eta_{j+1} - E[\eta_{j+1}|\eta_j]\big)^2\big].
\end{eqnarray*}
Then,
\begin{eqnarray*}
\frac{1}{N}M_N(t) &=& Y_N(t) - Y_N(0) - \frac{1}{N}\sum_{j=0}^{[Nt]-1} (E(D_{N}) (\eta_{j}+1)-\eta_j)\\\
&=&  Y_N(t) - Y_N(0) - \frac{1}{N}\sum_{j=0}^{[Nt]-1} ((1-\frac{4c}{N +2c} ) (\eta_{j}+1)-\eta_j)  \\
&=& Y_N(t) - Y_N(0) - \frac{1}{N}[Nt]  + \frac{4c}{N+2c}\sum_{j=0}^{[Nt]-1} Y_N\left(\frac{j}{N}\right) +   \frac{4c[Nt]}{N(N+2c)}
 \end{eqnarray*}
and
\begin{eqnarray}
\langle N^{-1}M_N(t)\rangle &=& \frac{1}{N^2} \sum_{j=0}^{[Nt]-1} E\big[(\eta_{j+1}-E(D_{N})(\eta_j+1))^2|\eta_j\big]\nonumber\\
&=& \frac{1}{N^2} \sum_{j=0}^{[Nt]-1} E\left[\left(\sum_{i=0}^{\eta_j} (\xi_{j,i}-E(D_{N}))\right)^2|\eta_j\right]\nonumber\\
&=& \frac{1}{N^2} \sum_{j=0}^{[Nt]-1} (\eta_j +1) \mbox{Var}(D_{N})\nonumber\\
& = &\frac{2}{N}\sum_{j=0}^{[Nt]-1}Y_N\left(\frac{j}{N}\right) + O\left(\frac{1}{N}\right).
\label{quad_var1}
\end{eqnarray}

Now suppose $Y_N(t)$ and $N^{-1}M_N(t)$ are tight in the sup topology, and $Y_N(t)\Rightarrow Z_t$, $N^{-1}M_N(t)\Rightarrow M(t)$ on subsequences.  Then,
$M(t)  =  Z_t-Z_0 -\int_{0}^{t} (1-4cZ_{s})ds$ and $\langle M(t)\rangle  =  2\int_0^t Z_s ds$.
Hence, by Levy's criterion for continuous martingales, we have that $Z_t$ is uniquely characterized by
$$Z_t \ = Z_{0} + \int_{0}^{t} (1-4cZ_{s})ds + \int_{0}^{t}\sqrt{2Z_{s}}dB_s$$

Similarly, for $[Nt]\geq [(1-\alpha) N]$, since now $\eta_{j+1}= \sum_{i=1}^{\eta_j}\xi_{j,i}$, the drift is not present, and we can write
\begin{eqnarray*}
\frac{1}{N}(M_N(t)- M_N(1-\alpha)) & = &  Y_N(t) - Y_N(1-\alpha) + \sum_{j = [(1-\alpha)N]}^{[Nt] - 1} (E(\eta_{j+1} | \eta_{j}) - \eta_{j})\\
 &=&  Y_N(t) - Y_N(1-\alpha) + \frac{4c}{N +c}\sum_{j=[N(1-\alpha)]}^{[Nt]-1} Y_N\left(\frac{j}{N}\right)
 \end{eqnarray*}
 and also
 \begin{eqnarray*}
      \frac{1}{N}(\langle M_N(t)\rangle -\langle M_N(1-\alpha)\rangle)
&=& \frac{1}{N^2} \sum_{j=[N(1-\alpha)]}^{[Nt]-1} (\eta_j +1) \mbox{Var}(D_{N})\\
 & =&
\frac{2}{N} \sum_{j=[N(1-\alpha)]}^{[Nt]-1} Y_N\left(\frac{j}{N}\right) + O\left(\frac{1}{N}\right). \nonumber\\
\end{eqnarray*}
Hence, as before, given tightness of $N^{-1}M_N(t)$, and subsequential convergences $Y_n(t)\Rightarrow Z(t)$ and $N^{-1}M_N(t) \Rightarrow M(t)$ on $[1-\alpha, 1]$ where $M(t) - M(1-\alpha) = Z_t-Z_{1-\alpha}$ and $\langle M(t) - M(1-\alpha)\rangle = 2\int_{1-\alpha}^t Z_s ds$, we conclude, for $t\in [1-\alpha,1]$, that
$$Z_t \ = Z_{1-\alpha} + \int_{1-\alpha}^{t }\sqrt{2Z_{s}}dB_s.$$
Consequently, it follows, putting the subsequential converges together, for $0 \leq t \leq 1$ that $Y_N(t)$ converges weakly to $Z_{t}$.
\vskip .1cm

{\it Tightness.}  We now argue tightness of $Y_N(t)$ and $N^{-1}M_N(t)$ on $[0,1-\alpha]$.  Tightness of $Y_N(t)$ follows from tightness of $N^{-1}M_N(t)$ in the sup topology which can be argued by a Kolmogorov-Centsov argument. First, for a general discrete time martingale $(M(l), \mathcal F_l)$ with difference $\d(l) = M(l)-M(l-1)$, we have that
\begin{eqnarray*}
E\left[ (M(l)-M(k))^4\right] &=& 6\sum_{j=k+1}^l E\left[ \d(j)^2(M(j-1)-M(k))^2\right]\\
&&\ \ \ + 4 \sum_{j=k+1}^l E\left[ \d(j)^3(M(j-1)-M(k))\right] + \sum_{j=k+1}^l E\left[\d^4_j\right]
\end{eqnarray*}
and by Jensen inequality,
\begin{eqnarray*}
&&E\left[ (M(l)-M(k))^4\right] \ \leq \ 6\sum_{j=k+1}^l E\left[ E[\d(j)^2|\mathcal F_{j-1}](M(j-1)-M(k))^2\right]\\
&&\ \ \ \ \ \ \  +\  4 \sum_{j=k+1}^l \left\{E\left[ E[\d(j)^4|\mathcal F_{j-1}]^{3/2}(M(j-1)-M(k))^2\right]\right\}^{1/2}\\
 &&\ \ \ \ \ \ \  +\  \sum_{j=k+1}^l E\left[E[\d(j)^{4}|\mathcal F_{j-1}]\right].
\end{eqnarray*}

Now, in our context, define the martingale, for $l\leq [(1-\alpha)N]$,
$$M(l) \ = \ \eta_l -\eta_0 -\sum_{i=0}^{l-1} \left (E[\eta_{i+1}|\eta_i] - \eta_i\right )$$
so that $M_N(t) = M([Nt])$, and also the stopping time
$$\theta_{y,N} \ = \ \inf\{l\geq 0:  \eta_l \geq Ny\}.$$
Compute, with respect to $M(l\wedge \theta_{y,N})$, that
\begin{eqnarray*}
\d(l) & = & M(l\wedge \theta_{y,N})-M(l-1\wedge \theta_{y,N}) \\
& = & \eta_{l\wedge \theta_{N,y}} -  E[\eta_{l\wedge \theta_{y,N}}|\eta_{l-1\wedge \theta_{y,N}}]
 \ = \ \sum_{i=0}^{\eta_{l-1\wedge\theta_{N,y}}} (\xi_{l\wedge \theta_{y,N},i}-E(D_{N})).\end{eqnarray*}
Hence,
we have
$$E\Big[ \Big(\sum_{i=0}^{\eta_{l-1\wedge \theta_{y,N}}} (\xi_{l\wedge\theta_{y,N},i}-E(D_{N}))\Big)^2\Big |\mathcal F_{l-1\wedge\theta_{y,N}}\Big] \ \leq\  {\rm Var}(\X_{N})  \eta_{l-1\wedge \theta_{y,N}}\ \leq\  c_{1}Ny$$
and
\begin{eqnarray*}
E\Big[\Big(\sum_{i=0}^{\eta_{l-1\wedge \theta_{y,N}}} (\xi_{l\wedge \theta_{y,N},i}-E(D_{N}) )\Big)^4|\mathcal F_{l-1\wedge\theta_{y,N}}\Big] &\leq& c_{2}E\X_{N}^4  \eta^2_{l-1\wedge \theta_{y,N}} + \eta_{l-1\wedge \theta_{y,N}}\\
&\leq& (c_{2}E\X^4_N + 1) ((Ny)^2 + Ny).\end{eqnarray*}
Also, noting the quadratic variation estimate (\ref{quad_var1}),
$$\frac{1}{N^2}E\left[\left(M(j-1\wedge \theta_{N,y}\right)-M\left(k\wedge\theta_{N,y})\right)^2\right] \ \leq\  c_{3}|j-k|\left( \frac{1}{N} + y\right).$$

Hence, we have, for some constant $c_4$ not depending on $N$ or $y$, that
\begin{eqnarray*}
\frac{1}{N^4}E [ (M([Nt]\wedge\theta_{y,N}) - M([Ns]\wedge \theta_{y,N}))^4] &\leq& c_{4} \max\{y^2,1\}(|t-s|^2\vee \frac{1}{N^{2}}).
\end{eqnarray*}
Then, by Theorem 12.3 Billingsley \cite{Billingsley}, $N^{-1}M([Nt]\wedge \theta_{y,N})$ is tight for any $y<\infty$.  Hence, $N^{-1}M_N(t)$ is tight in the sup topology on $[0,1-\alpha]$.

Tightness with respect to the interval $[1-\alpha,1]$, and consequently the whole interval $[0,1]$ follows similarly. \qed

\subsection{Asymmetric walks}
The situation is much different for asymmetric walks, in particular, the local times are 
of  order $O(1)$, and no scaling is required.
Given $p>q$, the walk starting from $x=[\alpha N]$ will exit to the right with probability tending to $1$ as $N\uparrow\infty$.  The sequence $\eta_j$ for $1\leq j\leq [(1-\alpha) N]$ is a branching process with mean offspring $E\X_N = q/p<1$ and immigration at each time of one individual.  The initial population is $\eta_1$ with the distribution of $D=\X_N$, a Geometric$(q)$ random variable minus $1$.  Hence, this sequence is a positive recurrent Markov chain, and $\eta_{[(1-\alpha)N]}$ converges to the stationary distribution $\pi$.

 On the other hand, the chain $\eta_j$ for $j\geq [(1-\alpha) N]$ is the usual Branching process with offspring distribution $\X_N$ (and no immigration).  Hence, it dies out in finite time.

The stationary distribution $\pi$ can be described by its probability generating function $\Psi(s) = \sum_{k\geq 0} \pi(k)s^k$.  Let $\phi(s)$ be the probability generating function of $D$.  Then, easy computations give that $\Psi(s)=\Psi(\phi(s))\phi(s)$.

Hence, since the distribution of $\eta_{[(1-\alpha)N]}$ converges to $\pi$, we can state a limit characterization in terms of a reversed process.
\begin{proposition}
  \label{asym_limprop}  Consider the asymmetric walk when $p>q$ starting from $[\alpha N]$.
  For any $M\geq 1$,
 the reversed process $\{\beta_k=\eta_{[(1-\alpha) N]-k}\}_{k= 0}^M$ converges in distribution to the reversed process starting from the stationary distribution $\pi$ of the chain.

 However, $\{\beta_{-k}= \eta_{[(1-\alpha)N] +k}\}_{k=0}^M$ converges to a Branching process with offspring distribution $D$ starting from $\pi$.
\end{proposition}

%
%
%

\section{Question 3: Periodicity}
\label{independent}

We now address the parity of various well-separated locations visited by the walk before exiting.  We remark different types of multiple point structures in other settings have been studied in Hamana \cite{Ham} and Pitt \cite{Pitt}.


Let $0<\alpha_1<\alpha_2<\cdots<\alpha_k<1$, and $e_i \in \{0,1\}$ for $1\leq i\leq k$.
\begin{proposition}
\label{thm1_question3}
 With respect to symmetric or weakly asymmetric walks, for $\alpha\in (0,1)$, we have
$$\lim_{N\rightarrow\infty} P_{[\alpha N]}\left( \cap_{i=1}^k \{G([\alpha_i N]) = e_i \ {\rm mod}_2\}| \max_{1\leq i\leq k}T_{[\alpha_i N]}<\tau_N\right) \ = \ \frac{1}{2^k}.$$
\end{proposition}
In other words, in the symmetric or weakly asymmetric cases, given that the locations are visited, the parities at $\{[\alpha_i N]\}_{i=1}^k$ converge to i.i.d. fair Bernoulli random variables.

But, with respect to asymmetric walks when $p>q$, starting from $[\alpha N]$, unless $\alpha<\beta$, $[\beta N]$ is not visited with probability tending to $1$.  So, it makes sense only to discuss parities of sites to the right of $[\alpha N]$.
\begin{proposition}
\label{thm2_question3}
With respect to asymmetric walks when $p>q$, suppose $0<\alpha<\alpha_1$.  Then,
$$\lim_{N\rightarrow\infty} P_{[\alpha N]}\left( \cap_{i=1}^k \{G([\alpha_i N]) = 1 \ {\rm mod}_2\}\right) \ = \ \frac{1}{(2-(p-q))^k}.$$
\end{proposition}
By the inclusion-exclusion principle, one concludes, in the asymmetric situation, the parities at $\{[\alpha_i N]\}_{i=1}^k$ converge
to i.i.d. Bernoulli random variables with success probability $(2-(p-q))^{-1}$.

We remark, with respect to the `stochastic locker' intepretation, one concludes that the expected proportion of lockers left 
closed is half or $(2-(p-q))^{-1}$ times the proportion of the range in the symmetric/weakly asymmetric, or asymmetric cases respectively.

\subsection{Proofs of Propositions \ref{thm1_question3} and \ref{thm2_question3}}
The proofs of the above propositions are similar.  We first derive the chance a single site is left open, and then later use this development in an induction scheme.
Let $T^r_y$ be the $r$th hitting time of $y$, and $\tilde{T}_y = \inf\{n\geq 1: X_n = y\}$ be the return time to $y$.
The event that site $y$ is left open, with various prescribed exits, is expressed as
\begin{eqnarray*}
&&\{G(y) = 1 \ {\rm mod}_2, T_N<T_0\}\\
&&\ \ \ \ \ \ \ \ \  =\ \cup_{k\geq 0} \{T^1_y<\tau_N\}\cap \{T^{2k+1}_y<\tau_N\}\cap\{T_N<T^{2k+2}_y\wedge T_0\}.
\end{eqnarray*}
Similarly,
\begin{eqnarray*}
&&\{G(y) = 1 \ {\rm mod}_2, T_0<T_N\}\\
&&\ \ \ \ \ \ \ \ \ = \  \cup_{k\geq 0} \{T^1_y<\tau_N\}\cap \{T^{2k+1}_y<\tau_N\}\cap\{T_0<T^{2k+2}_y\wedge T_N\}\\
&&\{G(y) = 1 \ {\rm mod}_2\} \ = \ \cup_{k\geq 0} \{T^1_y<\tau_N\}\cap \{T^{2k+1}_y<\tau_N\}\cap\{\tau_N<T^{2k+2}_y\}.
\end{eqnarray*}

Then,
\begin{eqnarray*}
P_x(G(y) = 1 \ {\rm mod}_2, T_N<T_0) &=& P_x(T_y<\tau_N) P_y(T_N<\tilde T_y)\sum_{l\geq 0} P_y(\tilde T_y<\tau_N)^{2l}\\
&=&\frac{P_x(T_y<\tau_N) P_y(T_N<\tilde T_y)}{1-(1-P_y(\tau_N<\tilde T_y))^2}\\
&=&\frac{P_x(T_y<\tau_N)}{2-P_y(\tau_N<\tilde T_y)}\frac{P_y(T_N<\tilde T_y)}{P_y(\tau_N<\tilde T_y)}.
\end{eqnarray*}
Also,
\begin{eqnarray*}
P_x(G(y)=1\ {\rm mod}_2, T_0<T_N) &=&\frac{P_x(T_y<\tau_N)}{2-P_y(\tau_N<T_y)}\frac{P_y(T_0<\tilde T_y)}{P_y(\tau_N<\tilde T_y)}\\
P_x(G(y)=1\ {\rm mod}_2) &=&\frac{P_x(T_y<\tau_N)}{2-P_y(\tau_N<\tilde T_y)}.
\end{eqnarray*}
In this last expression $P_x(T_y<\tau_N)$ is the probability $y$ is visited starting from $x$, and $(2-P_y(\tau_N<\tilde T_y))^{-1}$ is the factor specifying that $y$ is left open.  The quantity $P_y(\tau_N<\tilde T_y)$ can be viewed as an ``escape probability.''

Suppose now $x=[ \alpha N]$ and $y=[ \beta N]$.  In the symmetric case, we compute
$$P_x(T_y<\tau_N) = \left\{\begin{array}{rl}
\frac{N-x}{N-y} & {\rm \ for \ }y<x<N, \\
\frac{x}{y} &\ {\rm for \ }0<x<y\end{array}\right.$$
and
$$
P_y(\tilde T_y<\tau_N) \ = \ \frac{1}{2}P_{y-1}(T_y<T_0) +
\frac{1}{2}P_{y+1}(T_y<T_N) \ = \ 1-\frac{N}{2y(N-y)}.$$

In the (weakly) asymmetric case, we have
$$
P_x(T_y<\tau_N) = \left\{\begin{array}{rl}
\frac{s_N^x-s_N^N}{s_N^y-s_N^N} & \ {\rm for \ }x>y\\
\frac{1-s_N^x}{1-s_N^y}&\ {\rm for \ }x<y\end{array}\right.$$
and
\begin{eqnarray*} P_y(\tau_N<\tilde T_y) &=& q_NP_{y-1}(T_0<T_y) + p_NP_{y+1}(T_N<T_y) \\
&=&\frac{p_N(1-s_N)(1-s_N^N)}{(1-s_N^y)(1-s_N^{N-y})}.\end{eqnarray*}
Then,
$$P_y(\tau_N<\tilde T_y) \ \rightarrow\  \left\{\begin{array}{rl}
0 & \ {\rm for \ symmetric/weakly \ asymmetric \ walks }\\
p-q& \ {\rm for \ asymmetric \ walks}.
\end{array}\right.
$$

Putting these observations together, we have the following result.
\begin{proposition}
\label{oddprop}
Under symmetric or weakly asymmetric motion,
\begin{eqnarray*}
&&\lim_{N\uparrow\infty}P_x(G(y)=1\ {\rm mod}_2| T_y<T_N<T_0)\\
&&\ \ \ \ \ \ \ \ \  =\  \lim_{N\uparrow\infty}P_x(G(y)=1\ {\rm mod}_2| T_y<T_0< T_N) \  =\  \frac{1}{2},
\end{eqnarray*}
and hence $\lim_{N\uparrow\infty}P_x(G(y)=1\ {\rm mod}_2| T_y<\tau_N) = 1/2$.

However, under asymmetric motion, for $x\leq y$,
$$\lim_{N\uparrow\infty}P_x(G(y)=1\ {\rm mod}_2)  \ = \ \frac{1}{2-(p-q)}.$$
\end{proposition}


%
\vskip .1cm

{\it Proof of Proposition \ref{thm1_question3}.}
Let $G_n(y) = \sum_{l=0}^{n\wedge \tau_N} 1_y(X_l)$ be the number of visits to $y$ up to time $n\wedge \tau_N$.  First, we write
\begin{eqnarray}
&&P_{[\alpha N]}(\cap_{i=1}^k \{T_{[\alpha_i N]}<\tau_N\},\cap_{i=1}^k G([\alpha_i N])=e_i \ {\rm mod}_2)\nonumber\\
 &&\ \ = \ P_{[\alpha N]}(\cap_{i=1}^k \{G([\alpha_i N])=e_i \ {\rm mod}_2\}, T_{[\alpha_1 N]}<T_N<T_0)\nonumber\\
 &&\ \ \ \ \ \
+ P_{[\alpha N]}(\cap_{i=1}^k G([\alpha_i N])=e_i \ {\rm mod}_2, T_{[\alpha_k N]}<T_0<T_N).
\label{sec_3_eqn0}
\end{eqnarray}

We now concentrate on the first term on the right when $T_N<T_0$, as the argument is similar for the second term.
Since, on the set $T_N<T_0$, the walk must leave $[\alpha_1 N]$ never to return, and is also nearest-neighbor, write
\begin{eqnarray}
&&P_{[\alpha N]}(T_{[\alpha_1 N]}<T_N<T_0, \cap_{i=1}^k\{G([\alpha_i N])=e_i \  {\rm mod}_2\})\nonumber\\
&&\ \ \ \ = \ \sum_{\stackrel{z_1,\ldots, z_k}{z_1 = e_1 \ {\rm mod}_2}} P_{[\alpha N]}(T^{z_1}_{[\alpha_1 N]}<\tau_N, \cap_{i=2}^k \{G_{T^{z_1}_{[\alpha 1 N]}}([\alpha_i N]) = z_i\})\nonumber\\
&&\ \ \ \ \ \ \ \ \ \cdot \ P_{[\alpha_1 N]}( T_{[\alpha_2 N]}<\tilde T_{[\alpha_1 N]}\wedge T_N)\nonumber\\
&&\ \ \ \ \ \ \ \ \ \cdot \ P_{[\alpha_2 N]}(T_N<T_{[\alpha_1 N]}, \cap_{i=2}^k \{G([\alpha_i N])=e(z_i)\})
\label{sec_3_eqn}
\end{eqnarray}
where $e(z_i)=e_i$ or $1-e_i$ if $z_i$ is even or odd respectively.

In the last factor, which deals with the parities of $k-1$ points, $[\alpha_1 N]$ can be translated to $x=0$.  Treating the limit in Proposition \ref{oddprop} as a base step, we may conclude by induction, for fixed $e(z_i)$, that
$$\lim_{N\rightarrow \infty}P_{[\alpha_2 N]}(\cap_{i=2}^k \{G([\alpha_i N])=e(z_i)\}|T_N<T_{[\alpha_1 N]}) \ = \ {2}^{-(k-1)}.$$

Hence, by bounded convergence, we may replace the last factor of (\ref{sec_3_eqn}), by
$$2^{-(k-1)}P_{[\alpha_2N]}(T_N<T_{[\alpha_1N]}) + o(1)$$
as $N\uparrow\infty$.  Summing over $z_2,\ldots, z_k$, we have
\begin{eqnarray*}
&&P_{[\alpha N]}(T_{[\alpha_1 N]}<T_N<T_0, \cap_{i=1}^k\{G([\alpha_i N])=e_i\})\\
&& \ \ = \  \left[\sum_{z_1 = e_1 \ {\rm mod}_2}P_{[\alpha N]}(T^{z_1}_{[\alpha_1 N]}<\tau_N)\right ]\cdot P_{[\alpha_1 N]}(T_{[\alpha_2 N]}<\tilde T_{[\alpha_1 N]}\wedge T_N)\\
&&\ \ \ \ \ \ \ \ \ \ \ \ \  \cdot [2^{-(k-1)}P_{[\alpha_2N]}(T_N<T_{[\alpha_1N]}) +o(1)] \\
&& \ \ = \  \frac{1}{2^{k-1}}P_{[\alpha N]}(G([\alpha_1 N])= e_1 \ {\rm mod}_2, T_{[\alpha_1 N]}<T_N<T_0) +o(1).
\end{eqnarray*}

Therefore, noting Proposition \ref{oddprop},
$$\lim_{N\rightarrow\infty}
P_{[\alpha N]}(\cap_{i=1}^k\{G([\alpha_i N])=e_i\}| T_{[\alpha_1 N]}<T_N<T_0) \ = \ \frac{1}{2^k}.
$$
A similar expression is derived when the conditioning event is $T_{[\alpha_k N]}<T_0<T_N$, and so the limit in Proposition \ref{thm1_question3} is recovered. \qed

\vskip .1cm
{\it Proof of Proposition \ref{thm2_question3}.}  The proof is easier than that for Proposition \ref{thm1_question3}.
Since the probability of backtracking, $P_{[\gamma N]}(T_{[\beta N]}<\tau_N)$ is exponentially small in $N$ for $\beta<\gamma$, and noting Proposition \ref{oddprop}, we have
\begin{eqnarray*}
P_{[\alpha N]}(\cap_{i=1}^k \{G([\alpha_i N]) = 1 \ {\rm mod}_2\}) &=&
o(1)+\prod_{i=1}^k P_{[\alpha_i N]}(G([\alpha_i N]) = 1 \ {\rm mod}_2)\\
&\rightarrow & (2-(p-q))^{-k}
\end{eqnarray*}
\qed

\end{document}